# CONTROL OF GENERALIZED ERROR RATES IN MULTIPLE TESTING

By Joseph P. Romano[1] and Michael Wolf[2]

*Stanford University and University of Zurich*

Consider the problem of testing $s$ hypotheses simultaneously. The usual approach restricts attention to procedures that control the probability of even one false rejection, the familywise error rate (FWER). If $s$ is large, one might be willing to tolerate more than one false rejection, thereby increasing the ability of the procedure to correctly reject false null hypotheses. One possibility is to replace control of the FWER by control of the probability of $k$ or more false rejections, which is called the $k$-FWER. We derive both single-step and step-down procedures that control the $k$-FWER in finite samples or asymptotically, depending on the situation. We also consider the false discovery proportion (FDP) defined as the number of false rejections divided by the total number of rejections (and defined to be 0 if there are no rejections). The false discovery rate proposed by Benjamini and Hochberg [*J. Roy. Statist. Soc. Ser. B* **57** (1995) 289–300] controls $E(\text{FDP})$. Here, the goal is to construct methods which satisfy, for a given $\gamma$ and $\alpha$, $P\{\text{FDP} > \gamma\} \leq \alpha$, at least asymptotically. In contrast to the proposals of Lehmann and Romano [*Ann. Statist.* **33** (2005) 1138–1154], we construct methods that implicitly take into account the dependence structure of the individual test statistics in order to further increase the ability to detect false null hypotheses. This feature is also shared by related work of van der Laan, Dudoit and Pollard [*Stat. Appl. Genet. Mol. Biol.* **3** (2004) article 15], but our methodology is quite different. Like the work of Pollard and van der Laan [*Proc. 2003 International Multi-Conference in Computer Science and Engineering, METMBS'03 Conference* (2003) 3–9] and Dudoit, van der Laan and Pollard [*Stat. Appl. Genet. Mol. Biol.* **3** (2004) article 13], we employ resampling methods to achieve our

Received June 2005; revised November 2006.

[1]Supported by National Science Foundation Grant DMS-01-03926.

[2]Supported by the Spanish Ministry of Science and Technology and FEDER, Grant BMF2003-03324.

*AMS 2000 subject classifications.* Primary 62J15; secondary 62G10.

*Key words and phrases.* Bootstrap, false discovery proportion, false discovery rate, generalized familywise error rate, multiple testing, step-down procedure.

goals. Some simulations compare finite sample performance to currently available methods.

**1. Introduction.** The main goal of this paper is to show how computer-intensive methods can be used to construct asymptotically valid tests of multiple hypotheses under very weak conditions. In particular, we construct computationally feasible methods which provide control (at least asymptotically) of some generalized notions of the familywise error rate. However, the theory also applies to exact finite sample control in certain situations.

Consider the problem of testing hypotheses $H_1, \ldots, H_s$. A classical approach to dealing with the multiplicity problem is to restrict attention to procedures that control the probability of one or more false rejections, which is called the familywise error rate (FWER). For a given family, control of the FWER at (joint) level $\alpha$ requires that FWER $\leq \alpha$ for all possible distributions of the data considered in the model.

Of course, safeguards against false rejections are not the only concern of multiple testing procedures. Corresponding to the power of a single test, one must also consider the ability of a procedure to detect departures from the null hypotheses. When the number of tests, $s$, is large, such as in genomics studies, control of the FWER at conventional levels becomes so stringent that individual departures from the null hypotheses have little chance of being detected. For this reason, we shall consider alternatives to the FWER that control false rejections less severely in hopes of better power.

First, we shall consider the $k$-FWER, the probability of rejecting at least $k$ true null hypotheses. More formally, suppose data $X$ is available from some model $P \in \Omega$. A general hypothesis $H$ can be viewed as a subset $\omega$ of $\Omega$. For testing $H_i : P \in \omega_i$, $i = 1, \ldots, s$, let $I(P)$ denote the set of true null hypotheses when $P$ is the true probability distribution; that is, $i \in I(P)$ if and only if $P \in \omega_i$. Then, the $k$-FWER, which depends on $P$ is defined to be

$$(1) \qquad k\text{-FWER}_P = P\{\text{reject at least } k \text{ hypotheses } H_i : i \in I(P)\}.$$

Control of the $k$-FWER requires that $k$-FWER $\leq \alpha$ for all $P$; that is,

$$(2) \qquad\qquad k\text{-FWER}_P \leq \alpha \qquad \text{for all } P.$$

Evidently, the case $k = 1$ reduces to control of the usual FWER.

We will also consider control of the *false discovery proportion* (FDP), defined as the total number of false rejections divided by the total number of rejections (and equal to 0 if there are no rejections). Given a user specified value $\gamma \in [0, 1)$, the measure of error control we wish to control is $P\{\text{FDP} > \gamma\}$; thus, we wish to construct methods satisfying

$$(3) \qquad\qquad P\{\text{FDP} > \gamma\} \leq \alpha \qquad \text{for all } P.$$

We will derive methods where this holds (at least asymptotically). Evidently, control of the FDP with $\gamma = 0$ reduces to the usual FWER. Control of the false discovery rate (FDR) requires that $E(\mathrm{FDP}) \leq \gamma$.

Recently there have been a number of methods that control generalized error rates which are less stringent than the FWER. A prominent such technique is the FDR controlling method of [1]. Additional methods that control the FDR are given in [2] and [30]. Genovese and Wasserman [10] study asymptotic procedures that control the FDP (and the FDR) in the framework of a random effects mixture model. These ideas are extended in [20], where in the context of random fields, the number of null hypotheses is uncountable. Korn et al. [15] provide methods that control both the $k$-FWER and FDP; they provide some justification for their methods, but they are limited to a multivariate permutation model. Alternative methods of control of the $k$-FWER and FDP are given in [34]; they include both finite sample and asymptotic results. Like our work, their approach implicitly accounts for the dependence between the tests with the goal of improved ability to detect false hypotheses; comparisons between the methods will be made in Section 5. Building upon work for control of the FWER in [9, 22] and [28], we employ resampling to achieve our goals, which does not require the use of the subset pivotality condition of [35]. A further key ingredient is the use of the so-called $k$-max statistic, initially suggested in [9] in the construction of a single-step procedure. Our procedures here can be seen as step-down improvements over such single-step methods. A further new method is given in [33].

Some methods that control the $k$-FWER and FDP are now briefly reviewed. Suppose that $p$-values $\hat{p}_1, \ldots, \hat{p}_s$ are available for testing $H_1, \ldots, H_s$. For $\hat{p}_i$ to be a $p$-value, it is required that, for all $u \in [0, 1]$ and all $P \in \omega_i$, $P\{\hat{p}_i \leq u\} \leq u$. Then, for any fixed $k$, the procedure that rejects $H_i$ if $\hat{p}_i \leq k\alpha/s$ controls the $k$-FWER at level $\alpha$, and can be viewed as a generalization of the Bonferroni procedure which uses $k = 1$; see [17]. It is an example of a *single-step* procedure, meaning any null hypothesis is rejected if its corresponding $p$-value is less than or equal to a common cutoff value.

Improvements are possible by considering a class of *step-down* procedures, which we now describe. Order the $p$-values by $\hat{p}_{(1)} \leq \hat{p}_{(2)} \leq \cdots \leq \hat{p}_{(s)}$, and let $H_{(1)}, \ldots, H_{(s)}$ denote the corresponding hypotheses. Let

(4) $$\alpha_1 \leq \alpha_2 \leq \cdots \leq \alpha_s$$

be constants. If $\hat{p}_{(1)} > \alpha_1$, reject no null hypotheses. Otherwise, if

(5) $$\hat{p}_{(1)} \leq \alpha_1, \ldots, \hat{p}_{(r)} \leq \alpha_r,$$

reject hypotheses $H_{(1)}, \ldots, H_{(r)}$, where the largest $r$ satisfying (5) is used. The procedure of [13] uses $\alpha_j = \alpha/(s - j + 1)$ and controls the FWER at level $\alpha$. For general $k$, consider the following generalized Holm step-down procedure described in (5), where now we specifically set

$$(6) \qquad \alpha_j = \begin{cases} \dfrac{k\alpha}{s}, & j \le k, \\ \dfrac{k\alpha}{s+k-j}, & j > k. \end{cases}$$

Of course, the $\alpha_j$ depend on $s$ and $k$, but we suppress this dependence in the notation. Then the step-down method described in (5) with $\alpha_j$ given by (6) controls the $k$-FWER; that is, (2) holds; see [14] and [17].

Turning to FDP control, [17] reason as follows. To develop a step-down procedure satisfying (3), let $F$ denote the number of false rejections. At step $j$, having rejected $j-1$ hypotheses, we want to guarantee $F/j \le \gamma$, that is, $F \le \lfloor \gamma j \rfloor$, where $\lfloor x \rfloor$ is the greatest integer $\le x$. So, if $k = \lfloor \gamma j \rfloor + 1$, then $F \ge k$ should have probability no greater than $\alpha$; that is, we must control the number of false rejections to be $\le k$. Therefore, we use the step-down constant $\alpha_j$ with this choice of $k$ (which now depends on $j$); that is, we use the step-down

$$(7) \qquad \alpha_j = \frac{(\lfloor \gamma j \rfloor + 1)\alpha}{s + \lfloor \gamma j \rfloor + 1 - j}.$$

Under certain dependence assumptions on the $p$-values, this method satisfies (3). Some more conservative methods that hold under no dependence assumptions are also developed in [17, 25] and [26]. Typically, these generalized Holm type of methods assume a least favorable joint distribution for the $p$-values. In contrast, here we implicitly try to estimate the joint distribution of $p$-values with the hope of greater ability to detect false hypotheses.

In general, we suppose that rejection of $H_i$ is based on large values of a test statistic $T_{n,i}$ (with the subscript $n$ used for asymptotic purposes). If a $p$-value $\hat{p}_i$ is available for testing $H_i$, one can take $T_{n,i} = -\hat{p}_i$. Then we restrict attention to tests that reject an intersection hypothesis $H_K$ when the $k$th largest of the test statistics $\{T_{n,i} : i \in K\}$ is large. In some problems, [19] show that such stepwise procedures are optimal in a certain sense, in the case $k = 1$. Here, our primary goal is to show how computationally feasible step-down procedures can be constructed quite generally that control the $k$-FWER and FDP under minimal conditions.

In Section 2 we show that, if we estimate critical values that have a monotonicity property, then the basic problem of constructing a valid multiple test procedure that controls the $k$-FWER can essentially be reduced to the problem of sequentially constructing critical values for (at most order $s$) single tests that control the usual Type 1 error. In particular, if finite sample methods which offer control of the Type 1 error are available for each of the individual tests, then this will immediately translate into control of the $k$-FWER. Otherwise, we can apply bootstrap and subsampling methods to achieve asymptotic control, as described in Section 3. Results for control of the FDP are obtained in Section 4. Comparisons with the augmentation procedures of [34] are discussed in Section 5. In Section 6 we present a simulation study to examine the finite sample performance of various methods. The simulations demonstrate that our methods outperform or are at least competitive with currently available methods. All proofs are collected in the Appendix.

**2. Basic results for control of the $k$-FWER.** Suppose data $X$ is generated from some unknown probability distribution $P$. In anticipation of asymptotic results, we may write $X = X^{(n)}$, where $n$ typically refers to the sample size. A model assumes that $P$ belongs to a certain family of probability distributions $\Omega$, though we make no rigid requirements for $\Omega$; it may be a parametric, semiparametric or a nonparametric model.

Consider the problem of simultaneously testing a hypothesis $H_i$ against $H_i'$, for $i = 1, \dots, s$. Of course, a hypothesis $H_i$ can be viewed as a subset $\omega_i$ of $\Omega$, in which case the hypothesis $H_i$ is equivalent to $P \in \omega_i$ and $H_i'$ is equivalent to $P \notin \omega_i$. For any subset $K \subset \{1, \dots, s\}$, define $H_K = \bigcap_{i \in K} H_i$ to be the *intersection* hypothesis that $P \in \bigcap_{i \in K} \omega_i$. We also assume a test of the individual hypothesis $H_i$ is based on a test statistic $T_{n,i}$, with large values indicating evidence against $H_i$.

Some further notation is required. Suppose $\{y_i : i \in K\}$ is a collection of real numbers indexed by a finite set $K$ having $|K|$ elements. Then, for $k \leq |K|$, $k\text{-}\max(y_i : i \in K)$ is used to denote the $k$th largest value of the $y_i$ with $i \in K$. So, if the elements $y_i$, $i \in K$, are ordered as $y_{(1)} \leq \cdots \leq y_{(|K|)}$, then $k\text{-}\max(y_i : i \in K) = y_{(|K|-k+1)}$.

2.1. *Single-step control of the $k$-FWER.* Throughout this section, $k$ is fixed. First, we briefly discuss a single-step approach to control the $k$-FWER, since it serves as a building block for the more powerful step-down procedures considered later. For any subset $K \subset \{1, \dots, s\}$, let $c_{n,K}(\alpha, k, P)$ denote an $\alpha$-quantile of the distribution of $k\text{-}\max(T_{n,i} : i \in K)$ under $P$. Concretely,

$$(8) \qquad c_{n,K}(\alpha, k, P) = \inf\{x : P\{k\text{-}\max(T_{n,i} : i \in K) \leq x\} \geq \alpha\}.$$

(We use the subscript $n$ for asymptotic purposes, though the priority in this section is to study nonasymptotic results.) For testing the intersection hypothesis $H_K$ with $K \subset \{1, \dots, s\}$, it is only required to approximate a critical value for $P \in \bigcap_{i \in K} \omega_i$. Because there may be many such $P$, we define

$$(9) \qquad c_{n,K}(1-\alpha, k) = \sup\left\{ c_{n,K}(1-\alpha, k, P) : P \in \bigcap_{i \in K} \omega_i \right\}.$$

[In order to define $c_{n,K}(\alpha, k)$, we implicitly assume $\bigcap_{i=1}^{s} \omega_i$ is not empty.]

Consider the idealized test that rejects any $H_i$ for which $T_{n,i} > c_{n,I(P)}(1-\alpha, k, P)$. This is a single-step method in that each $T_{n,i}$ is compared with a common cutoff. However, this is an idealization because the critical value $c_{n,I(P)}(1-\alpha, k, P)$ is in general unknown. Such a fictional test clearly controls the $k$-FWER at level $\alpha$. Indeed, if $|I(P)| < k$, then there is nothing to prove; otherwise,

$$P\{k \text{ or more false rejections}\}$$

$$= P\{k\text{-}\max(T_{n,i} : i \in I(P)) > c_{n,I(P)}(1-\alpha, k, P)\} \leq \alpha,$$

with equality if the distribution of $k\text{-}\max(T_{n,i} : i \in I(P))$ is continuous under $P$. Unfortunately, the test is unavailable as the critical value is in general unknown. One possible approach is to replace $c_{n,I(P)}(1-\alpha, k, P)$ by $c_{n,I(P)}(1-\alpha, k)$, but this still depends on $P$ through $I(P)$. Since $I(P)$ is unknown, a conservative approach would be to assume all hypotheses are true and replace $c_{n,I(P)}(1-\alpha, k)$ by $c_{n,A}(1-\alpha, k)$, where $A = \{1, \ldots, s\}$.

Unfortunately, in nonparametric problems, the sup in (9) may be formidable or impossible to calculate, and may be way too conservative anyway. Instead, another possibility is to replace the critical value $c_{n,I(P)}(1-\alpha, k, P)$ by some estimate $\hat{c}_{n,I(P)}(1-\alpha, k)$, which is at least consistent or conservative. In general, suppose $\hat{c}_{n,K}(1-\alpha, k)$ represents an approximation or estimate of the $1-\alpha$ quantile of the distribution of $k\text{-}\max(T_{n,i} : i \in K)$, at least valid when $H_i$ is true for $i \in K$. Bootstrap and subsampling methods offer viable general approaches, and will be used later. Such a single-step approach using the $k$-max statistic was also discussed in [9]. (Rather than formalizing the required conditions for consistency right now, we will later give explicit conditions for more powerful step-down methods.) A single-step approach would then be to replace $K$ by $A = \{1, \ldots, s\}$.

EXAMPLE 2.1 (*Multivariate normal mean*). Suppose $(X_1, \ldots, X_s)$ is multivariate normal with unknown mean $\mu = (\mu_1, \ldots, \mu_s)$ and known covariance matrix $\Sigma$ having $(i, j)$ component $\sigma_{i,j}$. Consider testing $H_i : \mu_i \leq 0$ versus $\mu_i > 0$. Let $T_{n,i} = X_i/\sqrt{\sigma_{i,i}}$, since the test that rejects for large $X_i/\sqrt{\sigma_{i,i}}$ is *UMP* for testing $H_i$. For $|K| \geq k$, $c_{n,K}(1-\alpha, k)$ is the $1-\alpha$ quantile of the distribution of $k\text{-}\max(T_{n,i} : i \in K)$ when $\mu = 0$. A single-step approach would reject any $H_i$ for which $T_{n,i} > c_{n,A}(1-\alpha, k)$, where $A = \{1, \ldots, s\}$. Since $c_{n,A}(1-\alpha, k) \geq c_{n,I(P)}(1-\alpha, k) \geq c_{n,I(P)}(1-\alpha, k, P)$, this procedure clearly controls the $k$-FWER.

More generally, suppose $H_i$ specifies $\{P : \theta_i(P) \leq 0\}$ for some real-valued parameter $\theta_i$. Let $\hat{\theta}_{n,i}$ be an estimate of $\theta_i(P)$. Also, let $T_{n,i} = \tau_n \hat{\theta}_{n,i}$ for some nonnegative (nonrandom) sequence $\tau_n \to \infty$. The sequence $\tau_n$ is introduced for later asymptotic purposes so that a limiting distribution for $\tau_n[\hat{\theta}_{n,i} - \theta_i(P)]$ exists. In typical situations, $\tau_n = n^{1/2}$.

For $K \subset \{1, \ldots, s\}$ with $|K| \geq k$, let $L_{n,K}(k, P)$ denote the distribution under $P$ of $k$-$\max(\tau_n[\hat{\theta}_{n,i} - \theta_i(P)]: i \in K)$, with corresponding cumulative distribution function $L_{n,K}(x, k, P)$ and $\alpha$-quantile $b_{n,K}(\alpha, k, P) = \inf\{x: L_{n,K}(x, k, P) \geq \alpha\}$. By the definition of these quantiles and using $k = 1$,

$$(10) \qquad \left\{ (\theta_i: i \in A): \max_{i \in A} \tau_n[\hat{\theta}_{n,i} - \theta_i] \leq b_{n,A}(1 - \alpha, 1, P) \right\}$$

is an exact $1 - \alpha$ level joint confidence region for the subset of parameters $\{\theta_i(P): i \in A\}$. That is, the probability that the entire subset $\{\theta_i(P): i \in A\}$ will be contained in (10) is greater than or equal to $1 - \alpha$. By allowing $k \geq 1$, we can construct a "generalized" joint confidence region. More precisely, the probability that at least $|A| - k + 1$ elements of $\{\theta_i(P): i \in A\}$ will be contained in

$$(11) \qquad \left\{ (\theta_i: i \in A): k\text{-}\max_{i \in A} \tau_n[\hat{\theta}_{n,i} - \theta_i] \leq b_{n,A}(1 - \alpha, k, P) \right\}$$

is greater than or equal to $1 - \alpha$. In other words, the probability that $k$ or more elements of $\{\theta_i(P): i \in A\}$ will fall outside (11) is less than or equal to $\alpha$.

A value of 0 for $\theta_i(P)$ falls outside the region (11) if and only if $\tau_n \hat{\theta}_{n,i} > b_{n,A}(1 - \alpha, k, P)$. By the usual duality of confidence sets and hypothesis tests, this suggests the use of the critical value

$$(12) \qquad c_{n,A}(1 - \alpha, k) = b_{n,A}(1 - \alpha, k, P)$$

to control the $k$-FWER. The problem is that the critical value (12) is not feasible, since $P$ is unknown. Section 3 describe how approximate but feasible critical values can be obtained by the use of resampling methods. For example, the bootstrap replaces $P$ by an estimated distribution $\hat{Q}_n$, resulting in the critical value $\hat{c}_{n,A}(1 - \alpha, k) = b_{n,A}(1 - \alpha, k, \hat{Q}_n)$.

### 2.2. *Step-down methods that control the $k$-FWER.* Let

$$(13) \qquad T_{n,r_1} \geq T_{n,r_2} \geq \cdots \geq T_{n,r_s}$$

denote the observed ordered test statistics, and let $H_{r_1}, H_{r_2}, \ldots, H_{r_s}$ be the corresponding hypotheses. Step-down methods begin by first applying a single-step method, but then additional hypotheses may be rejected after this first stage by proceeding in a stepwise fashion, which we now describe. Begin by testing the joint null (intersection) hypothesis $H_{\{1,\ldots,s\}}$ that all hypotheses are true. This hypothesis is rejected if $T_{n,r_1}$ is deemed large, in which case $H_{r_1}$ is rejected. Here, the meaning of large is determined by some critical value $\hat{c}_{n,A}(1 - \alpha, k)$, which is designed to offer single-step control when testing the intersection hypothesis $H_A$ with $A = \{1, \ldots, s\}$. If it is not large, accept all hypotheses; otherwise, reject the hypothesis corresponding to the largest test statistic. Once a hypothesis is rejected, the next most significant hypothesis corresponding to the next largest test statistic is considered, and so on. At any stage, one tests appropriate intersection hypotheses $H_K$. Suppose that critical constants $\hat{c}_{n,K}(1 - \alpha, k)$ are available from our statistical tool chest, which we might contemplate for use as a single step procedure for testing $H_K$. The critical constants $\hat{c}_{n,K}(1 - \alpha, k)$ may be fixed or random, but the reader should have in mind that they each could be used as a test of $H_K$.

ALGORITHM 2.1 (*Generic step-down method for control of the k-FWER*).

1. Let $A_1 = \{1, \ldots, s\}$. If $\max(T_{n,i} : i \in A_1) \leq \hat{c}_{n,A_1}(1 - \alpha, k)$, then accept all hypotheses and stop; otherwise, reject any $H_i$ for which $T_{n,i} > \hat{c}_{n,A_1}(1 - \alpha, k)$ and continue.

2. Let $R_2$ be the indices $i$ of hypotheses $H_i$ previously rejected, and let $A_2$ be the indices of the remaining hypotheses. If $|R_2| < k$, then stop. Otherwise, let

$$\hat{d}_{n,A_2}(1 - \alpha, k) = \max_{I \subset R_2, |I| = k-1} \{\hat{c}_{n,K}(1 - \alpha, k) : K = A_2 \cup I\}.$$

Then, reject any $H_i$ with $i \in A_2$ satisfying $T_{n,i} > \hat{d}_{n,A_2}(1 - \alpha, k)$. If there are no further rejections, stop.

$$\vdots$$

j. Let $R_j$ be the indices $i$ of hypotheses $H_i$ previously rejected, and let $A_j$ be the indices of the remaining hypotheses. Let

$$\hat{d}_{n,A_j}(1 - \alpha, k) = \max_{I \subset R_j, |I| = k-1} \{\hat{c}_{n,K}(1 - \alpha, k) : K = A_j \cup I\}.$$

Then, reject any $H_i$ with $i \in A_j$ satisfying $T_{n,i} > \hat{d}_{n,A_j}(1 - \alpha, k)$. If there are no further rejections, stop.

$$\vdots$$

And so on.

Note that, in the case $k = 1$, once a hypothesis is removed, it no longer enters into the algorithm. However, for $k > 1$, the algorithm becomes slightly more complex. The reason is that, for control of the $k$-FWER, we must acknowledge that when we consider a set of hypotheses not previously rejected, we may have gotten to that stage by rejecting true null hypotheses, but hopefully at most $k-1$ of them. Since we do not know which of the hypotheses rejected thus far are true or false, we must maximize over subsets including some of those rejected, but at most $k-1$ among the previously rejected ones. Our main point will be that, if we can control the $k$-FWER at any stage of the algorithm, then the step-down test will control the $k$-FWER.

REMARK 2.1 (*Modified generic step-down method for control of the $k$-FWER*). One can modify the above algorithm or any method that controls the $k$-FWER as follows. If the method rejects at least $k-1$ hypotheses, no modification is applied; otherwise, reject the $k-1$ most significant hypotheses. This would not change control of the $k$-FWER. However, we do not generally promote this modification, because hypotheses can be rejected without compelling evidence (i.e., even if they have large unadjusted $p$-values).

In order to prove such an algorithm controls the $k$-FWER for suitable choice of critical values $\hat{c}_{n,K}(1-\alpha, k)$, we assume monotonicity of the estimated critical values; that is, for any $K \supset I(P)$,

$$(14) \qquad \hat{c}_{n,K}(1-\alpha, k) \geq \hat{c}_{n,I(P)}(1-\alpha, k).$$

Ideally, we would also like the following to hold: if $\hat{c}_{n,K}(1-\alpha, k)$ is used to test the intersection hypothesis $H_K$, then the chance of $k$ or more false rejections is bounded above by $\alpha$ when $K = I(P)$; that is,

$$(15) \qquad P\{k\text{-}\max(T_{n,i} : i \in I(P)) > \hat{c}_{n,I(P)}(1-\alpha, k)\} \leq \alpha.$$

Under the monotonicity assumption (14), we will show the basic inequality that $k$-FWER$_P$ is bounded above by left-hand side of (15). This will then show that, if we can construct monotone critical values such that each intersection test controls the $k$-FWER, then the step-down procedure controls the $k$-FWER. Thus, the construction of a step-down procedure is effectively reduced to construction of single tests, as long as the monotonicity assumption holds (and it always does for specific choices studied later).

THEOREM 2.1. *Let $P$ denote the true distribution generating the data. Consider Algorithm* 2.1 *with critical values $\hat{c}_{n,K}(1-\alpha, k)$ satisfying* (14).

(i) *Then*

$$(16) \qquad k\text{-}FWER_P \leq P\{k\text{-}\max(T_{n,i} : i \in I(P)) > \hat{c}_{n,I(P)}(1-\alpha, k)\}.$$

(ii) *Therefore, if the critical values also satisfy* (15), *then $k$-FWER$_P \leq \alpha$.*

The monotonicity assumption (14) cannot be removed, as shown in Example 2.1 of [28] in the case $k = 1$; an analogous construction works for general $k$. The general resampling constructions we describe later will inherently satisfy (14).

As a corollary, consider the nonrandom choice of critical values $\hat{c}_{n,K}(1 - \alpha, k) = c_{n,K}(1 - \alpha, k)$ defined in (9). Assume the following monotonicity assumption: for $K \supset I(P)$,

$$(17) \qquad c_{n,K}(1 - \alpha, k) \geq c_{n,I(P)}(1 - \alpha, k).$$

The condition (17) can be expected to hold in many situations because the left-hand side is based on computing the $1 - \alpha$ quantile of the $k$th largest of $|K|$ variables, while the right-hand side is based on the $k$th largest of $|I(P)| \leq |K|$ variables (though one must be careful and realize that the quantiles are computed under possibly different $P$, which is why some condition is required).

COROLLARY 2.1.  *Let $P$ denote the true distribution generating the data. Assume $\bigcap_{i=1}^{s} \omega_i$ is not empty.*

(i) *Consider Algorithm 2.1 with $\hat{c}_{n,K}(1 - \alpha, k) = c_{n,K}(1 - \alpha, k)$ and assume* (17). *Then $k$-$FWER_P \leq \alpha$.*

(ii) *Control persists if in Algorithm 2.1 the critical constants $\hat{c}_{n,K}(1 - \alpha, k)$ are replaced by $d_{n,K}(1 - \alpha, k)$ which satisfy $d_{n,K}(1 - \alpha, k) \geq c_{n,K}(1 - \alpha, k)$.*

(iii) *Moreover, the condition* (17) *may then be removed if the $d_{n,K}(1 - \alpha, k)$ satisfy $d_{n,K}(1 - \alpha, k) \geq d_{n,I(P)}(1 - \alpha, k)$ for any $K \supset I(P)$.*

EXAMPLE 2.2 (*Multivariate normal mean, continuation of Example 2.1*). Recall the setup of Example 2.1 with $T_{n,i} = X_i / \sqrt{\sigma_{i,i}}$. To apply Corollary 2.1, assume that $|I(P)| \geq k$ or there is nothing to prove. Let $c_{n,K}(1 - \alpha, k)$ be the $1 - \alpha$ quantile of the distribution of $k$-$\max(T_{n,i} : i \in K)$ when $\mu = 0$. Since $k$-$\max(T_{n,i} : i \in I) \leq k$-$\max(T_{n,i} : i \in K)$ whenever $I \subset K$, (17) is satisfied. Moreover, the resulting procedure rejects at least as many hypotheses as the generalized Holm procedure, as it accounts for the dependence of the test statistics.

The previous example is parametric in nature. However, we will see that a valid step-down approach can apply to nonparametric problems. Our main goal will be to apply resampling methods that can account for the dependence structure of the test statistics. We also observe that Theorem 2.1 applies to certain semiparametric problems where permutation and randomization tests apply. This was accomplished in the case $k = 1$ by [28], but the argument generalizes given Theorem 2.1. In fact, the result in [15] for $k$-FWER control in a specialized multivariate permutation setup is a special case of our results.

However, we first observe the fact that the generalized Holm procedure described by (5) with critical values given by (6) controls the $k$-FWER. This follows from Theorem 2.1 and the fact that, when testing $|K|$ hypotheses, the single-step procedure that rejects any $H_i$ whose corresponding $p$-value is $\leq k\alpha/|K|$ controls the $k$-FWER; see Theorem 2.1(i) of [17]. Note the critical values $k\alpha/|K|$ are monotone in $|K|$.

Outside some parametric models, application of the generic step-down method can be computationally intensive, so we will also consider the following more streamlined algorithm. The basic idea is that at any stage, when testing whether or not to include further rejections, we need only look at the hypotheses not previously rejected together with the $k-1$ hypotheses that are least significant among those previously rejected. So, we avoid maximizing over all subsets of size $k-1$ of previously rejected hypotheses and just look at the most "recent" $k-1$ rejections. The arguments for such a procedure will be asymptotic.

ALGORITHM 2.2 (*Streamlined step-down method for control of the $k$-FWER*). The algorithm is analogous to Algorithm 2.1. The only difference is that in any step $j > 1$ the critical value

$$\hat{d}_{n,A_j}(1-\alpha,k) = \max_{I \subset R_j, |I|=k-1}\{\hat{c}_{n,K}(1-\alpha,k) : K = A_j \cup I\}$$

is replaced by the critical value

$$\tilde{d}_{n,A_j}(1-\alpha,k) = \hat{c}_{n,K}(1-\alpha,k),$$

$$\text{where } K = \{r_{(|R_j|-k+2)}, r_{(|R_j|-k+1)}, \ldots, r_{(s)}\}.$$

**3. Asymptotic results on $k$-FWER control.** The main goal of this section is to show how Theorem 2.1 can be used to construct step-down procedures that asymptotically control the $k$-FWER under very weak assumptions. The use of resampling techniques will be a key ingredient. The methods constructed will be based on Algorithm 2.1, and so potentially many tests are constructed in a stepwise fashion. However, a key feature is that the methods will only require *one* set of resamples for all of the tests, whether they are bootstrap samples or subsamples.

In order to accomplish this, we will consider resampling schemes that do *not* obey the null hypothesis constraints. Such schemes have been suggested previously by [9] and [22], and have the benefit of avoiding the subset pivotality condition of [35]. Hypothesis test constructions that do obey the constraints imposed by the null hypothesis, as discussed in [4] and [24], are based on the idea that the critical value should be obtained under the null hypothesis and so the resampling scheme should reflect the constraints of the null hypothesis. This idea is even advocated as a principle in [12], and it is enforced throughout [35]. While appealing, it is by no means the only approach toward inference in hypothesis testing. In some problems, the *subset pivotality* condition of [35] holds, and so the same null distribution can be used at each step. However, this condition does not hold in general; for instance, see Example 4.1 of [28]. To obtain a more general construction, we exploit the well-known explicit duality between tests and confidence intervals; so, if one can construct good or valid confidence intervals, then one can construct good or valid tests, and conversely. The same holds for simultaneous confidence sets and multiple tests.

We shall consider two concrete applications of Theorem 2.1, the first based on the bootstrap and the second based on subsampling. The symbols $\overset{L}{\to}$ and $\overset{P}{\to}$ will denote convergence in law (or distribution) and convergence in probability, respectively.

3.1. *A bootstrap construction.* We now apply Theorem 2.1 to develop an asymptotically valid approach based on the bootstrap, but specializing to the case where $H_i$ is concerned with a test of a parameter. Suppose hypothesis $H_i$ is specified by $\{P : \theta_i(P) \leq 0\}$ for some real-valued parameter $\theta_i$. Implicitly, the alternatives are one-sided, but the two-sided case can be similarly handled. Suppose $\hat{\theta}_{n,i}$ is an estimate of $\theta_i$. Also, let $T_{n,i} = \tau_n \hat{\theta}_{n,i}$ for some nonnegative (nonrandom) sequence $\tau_n \to \infty$. The sequence $\tau_n$ is introduced for asymptotic purposes so that a limiting distribution for $\tau_n[\hat{\theta}_{n,i} - \theta_i(P)]$ exists. In typical situations, $\tau_n = n^{1/2}$.

The bootstrap method relies on its ability to approximate the joint distribution of $\{\tau_n[\hat{\theta}_{n,i} - \theta_i(P)] : i \in K\}$, which we denote by $J_{n,K}(P)$. For $K \subset \{1, \ldots, s\}$ with $|K| \geq k$, let $L_{n,K}(k, P)$ denote the distribution under $P$ of $k\text{-}\max(\tau_n[\hat{\theta}_{n,i} - \theta_i(P)] : i \in K)$, with corresponding c.d.f. $L_{n,K}(x, k, P)$ and $\alpha$-quantile $b_{n,K}(\alpha, k, P) = \inf\{x : L_{n,K}(x, k, P) \geq \alpha\}$.

We will assume the normalized estimates satisfy the following.

ASSUMPTION B1.

(i)  $J_{n,\{1,\ldots,s\}}(P) \overset{L}{\to} J_{\{1,\ldots,s\}}(P)$, a nondegenerate limit law.
(ii) $L_{I(P)}(\cdot, k, P)$ is continuous and strictly increasing on its support.

Part (i) implies that, for every $K \subset I(P)$, $L_{n,K}(k, P)$ has a limiting distribution $L_K(k, P)$. Indeed, the $k$-max function is a continuous function and the continuous mapping theorem applies; see Lemma A.1. Part (ii) makes an additional mild assumption on the limit law $L_{I(P)}(k, P)$. In particular, under Assumption B1, it follows that

$$
(18) \qquad b_{n, I(P)}(1 - \alpha, k, P) \to b_{I(P)}(1 - \alpha, k, P),
$$

where $b_{I(P)}(\alpha, k, P)$ is the $\alpha$-quantile of the limiting distribution $L_{I(P)}(k, P)$.

Let $\hat{Q}_n$ be some unrestricted estimate of $P$, that is, $\hat{Q}_n$ does not obey the null hypothesis constraints. For i.i.d. data, in the absence of a parametric model for $P$, $\hat{Q}_n$ is typically taken to be the empirical distribution of the observed data, or possibly a smoothed version (i.e., nonparametric bootstrap); on the other hand, if a parametric model for $P$ is assumed, then $\hat{Q}_n$ should be based on this model (i.e., parametric bootstrap); see [7]. For time series or data-dependent situations, bootstrap methods that can capture the underlying dependence structure should be employed, such as block bootstraps, sieve bootstraps or Markov bootstraps; see [16]. Then a nominal $1 - \alpha$ level bootstrap joint confidence region for the subset of parameters $\{\theta_i(P) : i \in K\}$ is given by

$$
\begin{aligned}
(19) \qquad & \{(\theta_i : i \in K) : \max(\tau_n[\hat{\theta}_{n,i} - \theta_i] : i \in K) \le b_{n,K}(1 - \alpha, 1, \hat{Q}_n)\} \\
& = \{(\theta_i : i \in K) : \theta_i \ge \hat{\theta}_{n,i} - \tau_n^{-1} b_{n,K}(1 - \alpha, 1, \hat{Q}_n)\}.
\end{aligned}
$$

So a value of 0 for $\theta_i(P)$ falls outside the region if and only if $\tau_n \hat{\theta}_{n,i} > b_{n,K}(1 - \alpha, 1, \hat{Q}_n)$. By the usual duality of confidence sets and hypothesis tests, this suggests the use of the critical value

$$
(20) \qquad \hat{c}_{n,K}(1 - \alpha, 1) = b_{n,K}(1 - \alpha, 1, \hat{Q}_n),
$$

to control the familywise error rate (i.e., the $k$-FWER with $k = 1$). Since here we require control of the $k$-FWER, we merely replace the max in (19) with the $k$-max and $b_{n,K}(1 - \alpha, 1, \hat{Q}_n)$ with $b_{n,K}(1 - \alpha, k, \hat{Q}_n)$. Such a generalized joint confidence region should asymptotically contain all true parameter values except for possibly at most $k - 1$ of them, with probability (asymptotically) at least $1 - \alpha$. Thus, the bootstrap critical value we use will be

$$
(21) \qquad \hat{c}_{n,K}(1 - \alpha, k) = b_{n,K}(1 - \alpha, k, \hat{Q}_n).
$$

Note that, regardless of asymptotic behavior, the monotonicity assumption (14) is always satisfied for the choice (21). Indeed, for any $Q$ and if $I \subset K$, $b_{n,I}(1 - \alpha, k, Q)$ is the $1 - \alpha$ quantile under $Q$ of the $k$-max of $|I|$ variables, while $b_{n,K}(1 - \alpha, k, Q)$ is the $1 - \alpha$ quantile of the $k$-max of these same $|I|$ variables together with $|K| - |I|$ additional variables. This simple observation together with Theorem 2.1 immediately yields:

COROLLARY 3.1. *Under the setup and notation of this subsection, consider Algorithm* 2.1 *with critical values given by* (21). *Then*

$$(22) \quad k\text{-FWER}_P \le P\{k\text{-}\max(T_{n,i} : i \in I(P)) > b_{n,I(P)}(1 - \alpha, k, \hat{Q}_n)\}.$$

Therefore, in order to conclude $\limsup_n k\text{-FWER}_P \le \alpha$, it is now only necessary to study the asymptotic behavior of $b_{n,I(P)}(1 - \alpha, k, \hat{Q}_n)$. For this, we further assume the usual conditions for bootstrap consistency when testing the *single* hypothesis that $\theta_i(P) \le 0$ for all $i \in I(P)$; that is, we assume the bootstrap consistently estimates the joint distribution of $\tau_n[\hat{\theta}_{n,i} - \theta_i(P)]$ for $i \in I(P)$. Specifically, consider the following (more general) assumption.

ASSUMPTION B2. *For any metric $\rho$ metrizing weak convergence on $\mathbb{R}^{|\{1,\dots,s\}|}$,*

$$\rho(J_{n,\{1,\dots,s\}}(P), J_{n,\{1,\dots,s\}}(\hat{Q}_n)) \xrightarrow{P} 0.$$

The Assumptions B1 and B2 are quite standard in the bootstrap literature, and readily hold for general classes of statistics, such as estimators which are smooth functions of means, $U$-statistics, $L$-statistics, estimators which are differentiable functions of the empirical process, and so forth; see [11, 31] and Chapter 1 of [21]. Thus, our results apply to a wide range of problems. Under these assumptions, the following theorem proves asymptotic control of the $k$-FWER of our bootstrap method.

THEOREM 3.1. *Fix $P$ satisfying Assumption* B1. *Let $\hat{Q}_n$ be an estimate of $P$ satisfying Assumption* B2. *Consider the method of Algorithm* 2.1 *with $\hat{c}_{n,K}(1 - \alpha, k)$ given by $b_{n,K}(1 - \alpha, k, \hat{Q}_n)$.*

(i) *Then $\limsup_n k\text{-FWER}_P \le \alpha$.*

(ii) *If $P$ is such that $i \notin I(P)$, that is, $H_i$ is false and $\theta_i(P) > 0$, then the probability that the step-down method rejects $H_i$ tends to 1.*

REMARK 3.1. Typically, one would like to choose test statistics that lead to procedures that are balanced in the sense that all tests have about the same power and contribute equally to error control, as argued by [5, 23] and [32]. Achieving balance is best handled by appropriate choice of test statistics. For example, using $p$-values as the basic statistics will lead to better balance. Quite generally, Beran's prepivoting transformation can lead to balance; see [5] and [6]. Alternatively, balance can sometimes be achieved by Studentization.

We now briefly consider the two-sided case. Suppose $H_i$ specifies $\theta_i(P) = 0$ against the alternative $\theta_i(P) \ne 0$. Let $L'_{n,K}(k, P)$ denote the distribution under $P$ of $k\text{-}\max(\tau_n|\hat{\theta}_{n,i} - \theta_i(P)| : i \in K)$ with corresponding distribution function $L'_{n,K}(x, k, P)$ and $\alpha$-quantile $b'_{n,K}(\alpha, k, P) = \inf\{x : L'_{n,K}(x, k, P) \geq \alpha\}$. Accordingly, $L'_K(k, P)$ denotes the limiting distribution of $L'_{n,K}(k, P)$. Finally, let $T'_{n,i} = \tau_n|\hat{\theta}_{n,i}|$. The following theorem extends Theorem 3.1 to the two-sided case.

THEOREM 3.2. *Fix $P$ satisfying Assumption B1, but with $L_{I(P)}(k, P)$ in Assumption B1(ii) replaced by $L'_{I(P)}(k, P)$. Let $\hat{Q}_n$ be an estimate of $P$ satisfying Assumption B2. Apply Algorithm 2.1 using the test statistics $T'_{n,i}$ and with $\hat{c}_{n,K}(1 - \alpha, k)$ given by $b'_{n,K}(1 - \alpha, k, \hat{Q}_n)$.*

(i) *Then $\limsup_n k\text{-}FWER_P \leq \alpha$.*

(ii) *If $P$ is such that $i \notin I(P)$, that is, $H_i$ is false and $\theta_i(P) \neq 0$, then the probability that the step-down method rejects $H_i$ tends to 1.*

(iii) *Moreover, if the above algorithm rejects $H_i$ and it is declared that $\theta_i > 0$ when $\hat{\theta}_{n,i} > 0$, the probability of making a Type 3 error [i.e., of declaring $\theta_i(P)$ positive when it is negative or declaring it negative when it is positive] tends to 0.*

So far, the bootstrap construction has been based on Algorithm 2.1. The following theorem shows that asymptotic control of the $k$-FWER is also achieved by the computationally less expensive streamlined Algorithm 2.2. For brevity we only focus on the one-sided case, that is, the setting of Theorem 3.1; the two-sided case is similar.

THEOREM 3.3. *Fix $P$ satisfying Assumption B1. Let $\hat{Q}_n$ be an estimate of $P$ satisfying Assumption B2. Consider the step-down method in Algorithm 2.2 with $\hat{c}_{n,K}(1 - \alpha, k)$ replaced by $b_{n,K}(1 - \alpha, k, \hat{Q}_n)$. Then the conclusions of Theorem 3.1 continue to hold.*

REMARK 3.2. The proofs of both Theorems 3.1 and 3.3 rely on asymptotic arguments. Nevertheless, some important differences should be pointed out. First, the method based on Algorithm 2.1 is more conservative than the one based on the streamlined Algorithm 2.2: the latter will reject all the hypotheses rejected by the former and potentially some further ones.

Second, if instead of the estimated critical values $b_{n,K}(1 - \alpha, k, \hat{Q}_n)$ the exact critical values $b_{n,K}(1 - \alpha, k, P)$ could be used in place of $\hat{c}_{n,K}(1 - \alpha, k)$, then Algorithm 2.1 would provide finite sample control of the $k$-FWER while Algorithm 2.2 would not.

Third, the bootstrap construction based on Algorithm 2.1 provides asymptotic control of the $k$-FWER in the case of contiguous alternatives while the construction based on Algorithm 2.2 may not. (An introduction to contiguity is given in Section 12.3 of [18].)

REMARK 3.3 (*Operative method*). The previous remark provides some motivation to base the bootstrap construction on the more conservative generic Algorithm 2.1. On the other hand, its computational burden can be very high. To compute the critical value $\hat{d}_{n,A_j}(1-\alpha, k)$ in the $j$th step, one has to evaluate $N_j = \binom{R_j}{k-1}$ quantiles $\hat{c}_{n,K}(1-\alpha, k)$ in order to then take the largest one of those. Depending on $R_j$ and $k$, this number $N_j$ may be very large. Therefore, we now suggest an operative method that retains some of the desirable properties of Algorithm 2.1 while remaining always computationally feasible. The suggestion is as follows. Pick a user specified number $N_{\max}$, say $N_{\max} = 50$, and let $M$ be the largest integer for which $\binom{M}{k-1} \le N_{\max}$. In step $j$ of Algorithm 2.1, the critical value is then computed as

$$\hat{d}_{n,A_j}(1-\alpha, k) = \max_{I \subset \{r_{\max\{1, |R_j|-M+1\}}, \ldots, r_{|R_j|}\}, |I|=k-1} \{\hat{c}_{n,K}(1-\alpha, k) : K = A_j \cup I\}.$$

That is, we maximize over subsets $I$ not necessarily of the entire index set $R_j$ of previously rejected hypotheses, but only of the index set corresponding to the $M$ least significant hypotheses rejected so far. (Of course, when $M \ge |R_j|$, we maximize over all subsets $I$ of $R_j$ of size $k - 1$.) The philosophy of this operative method is to be as close as possible to the generic Algorithm 2.1, given the limitation to the computational burden expressed by $N_{\max}$. Finally, note that the streamlined algorithm is a special case of the operative method when $N_{\max} = 1$ is chosen, resulting in $M = k - 1$.

**3.2. A general subsampling construction.** In this subsection, we present an alternative construction of critical values in our step-down procedure by using subsampling. Unlike the previous subsection, we do not assume $H_i$ is concerned with the test of a parameter $\theta_i$; the approach here is quite general and will hold under weaker asymptotic conditions as well. For any $K \subset \{1, \ldots, s\}$, let $G_{n,K}(P)$ be the joint distribution of the statistics $T_{n,i}$, $i \in K$, under $P$, with corresponding joint c.d.f. $G_{n,K}(x, P)$, $x \in \mathbb{R}^{|K|}$. Also, let $H_{n,K}(k, P)$ denote the distribution of $k\text{-}\max(T_{n,i} : i \in K)$ under $P$. As in Section 2.1, let $c_{n,K}(1-\alpha, k, P)$ denote a $1 - \alpha$ quantile of $H_{n,K}(k, P)$.

We will make the following general assumption.

ASSUMPTION S. Under $P$, the joint distribution of the test statistics $T_{n,i}$, $i \in I(P)$, has a limiting distribution; that is,

$$(23) \qquad\qquad G_{n,I(P)}(P) \xrightarrow{L} G_{I(P)}(P).$$

This implies that, under $P$, $k\text{-}\max(T_{n,i} : i \in I(P))$ has a limiting distribution, say $H_{I(P)}(k, P)$, with limiting c.d.f. $H_{I(P)}(x, k, P)$. Let $c_{I(P)}(\alpha, k, P)$ denote an $\alpha$-quantile of $H_{I(P)}(k, P)$; that is,

$$c_{I(P)}(\alpha, k, P) = \inf\{x : H_{I(P)}(x, k, P) \ge \alpha\}.$$

We will assume further that $H_{I(P)}(x, k, P)$ is continuous and strictly increasing at $x = c_{I(P)}(1 - \alpha, k, P)$.

Note that the above continuity condition is satisfied if the $|I(P)|$ univariate marginal distributions of $G_{I(P)}(P)$ are continuous; see Lemma A.1. Also, the strictly increasing assumption can be removed; see Remark 1.2.1 of [21].

We now detail the general subsampling construction. To this end, assume that we have available an i.i.d. sample $X_1, \ldots, X_n$ from $P$, and $T_{n,i} = T_{n,i}(X_1, \ldots, X_n)$ is the test statistic we wish to use for testing $H_i$. To describe the test construction, fix a positive integer $b < n$ and let $Y_1, \ldots, Y_{N_n}$ be equal to the $N_n := \binom{n}{b}$ subsets of $\{X_1, \ldots, X_n\}$, ordered in any fashion. Let $T_{b,i}^{(a)}$ be equal to the statistic $T_{b,i}$ evaluated at the data set $Y_a$, for $a = 1, \ldots, N_n$. Then, for any subset $K \subset \{1, \ldots, s\}$, the joint distribution of $(T_{n,i} : i \in K)$ can be approximated by the empirical distribution of the $N_n$ values $\{T_{b,i}^{(a)} : i \in K\}$. In other words, for $x \in \mathbb{R}^s$, the true joint c.d.f. of the test statistics evaluated at $x$,

$$G_{n,\{1,\ldots,s\}}(x, P) = P\{T_{n,1} \leq x_1, \ldots, T_{n,s} \leq x_s\},$$

is estimated by the subsampling distribution

$$(24) \qquad \hat{G}_{n,\{1,\ldots,s\}}(x) = \frac{1}{N_n} \sum_a I\{T_{b,1}^{(a)} \leq x_1, \ldots, T_{b,s}^{(a)} \leq x_s\}.$$

Note that the marginal distribution of any subset $K \subset \{1, \ldots, s\}$, $G_{n,K}(P)$, is then approximated by the marginal distribution induced by (24) on that subset of variables. So, $\hat{G}_{n,K}$ refers to the empirical distribution of the values $\{T_{n,i}^{(a)} : i \in K\}$. (In essence, one only has to estimate one joint sampling distribution for all the test statistics because this then induces that of any subset, even though we are not assuming anything like subset pivotality.)

Similarly, the estimate of the whole joint distribution of test statistics induces an estimate for the distribution of the maximum or $k$th largest of test statistics. Specifically, $H_{n,K}(k, P)$ is estimated by the empirical distribution $\hat{H}_{n,K}(x, k)$ of the values $k$-$\max(T_{n,i}^{(a)} : i \in K)$; that is,

$$\hat{H}_{n,K}(x, k) = \frac{1}{N_n} \sum_a I\{k\text{-}\max(T_{b,i}^{(a)} : i \in K) \leq x\}.$$

Also, let

$$(25) \qquad \hat{c}_{n,K}(1 - \alpha, k) = \inf\{x : \hat{H}_{n,K}(x, k) \geq 1 - \alpha\}$$

denote the estimated $1 - \alpha$ quantile of the $k$-max of test statistics $T_{n,i}$ with $i \in K$.

Note the monotonicity of the critical values: for $I \subset K$

$$(26) \qquad \hat{c}_{n,K}(1-\alpha, k) \geq \hat{c}_{n,I}(1-\alpha, k).$$

This simple observation together with Theorem 2.1 immediately yields:

COROLLARY 3.2. *Under the setup and notation of this subsection, consider Algorithm 2.1 with critical values given by (25). Then*

$$(27) \qquad k\text{-}FWER_P \leq P\{k\text{-}\max(T_{n,i} : i \in I(P)) > \hat{c}_{n,I(P)}(1-\alpha, k)\}.$$

The following result proves consistency and $k$-FWER control of our step-down algorithm based on these subsample estimates of critical values. Note, in particular, that Assumption B2 is not needed here at all, a reflection of the fact that the bootstrap requires much stronger (local uniform convergence) assumptions for consistency; see [21].

THEOREM 3.4. *Suppose Assumption S holds. Let $b/n \to 0$, $\tau_b/\tau_n \to 0$ and $b \to \infty$.*

   (i) *The subsampling approximation satisfies $\rho(\hat{G}_{n,I(P)}, G_{n,I(P)}(P)) \xrightarrow{P} 0$ for any metric $\rho$ metrizing weak convergence on $\mathbb{R}^{|I(P)|}$.*

   (ii) *The subsampling critical values satisfy $\hat{c}_{n,I(P)}(1-\alpha, k) \xrightarrow{P} c_{I(P)}(1-\alpha, k)$.*

   (iii) *Therefore, using Algorithm 2.1 with $\hat{c}_{n,K}(1-\alpha, k)$ given by (25) results in $\limsup_n k\text{-}FWER_P \leq \alpha$.*

The above approach can be extended to dependent data; see [21].

**4. Asymptotic results on FDP control.** In some applications, one might be willing to tolerate a larger number of false rejections in case the total number of rejections is large. In other words, one might be willing to tolerate a certain (small) fraction of false rejections out of the total rejections. This leads to control based on the *false discovery proportion* (FDP). Let $F$ be the number of false rejections made by a multiple testing procedure and let $R$ be the total number of rejections. Then the FDP is defined as

$$\text{FDP} = \begin{cases} \dfrac{F}{R}, & \text{if } R > 0, \\ 0, & \text{if } R = 0. \end{cases}$$

A multiple testing procedure is said to control the FDP at level $\alpha$ if, for the given sample size $n$, $P\{\text{FDP} > \gamma\} \leq \alpha$, for all $P$. A multiple testing procedure is said to asymptotically control the FDP at level $\alpha$, if $\limsup_n P\{\text{FDP} > \gamma\} \leq \alpha$, for all $P$. Our focus will be on procedures that provide asymptotic control. Notice that a procedure satisfying $P\{\text{FDP} > \gamma\} \leq 0.5$ guarantees that the median of the FDP is $\leq \gamma$. The main goal of this section is to construct a method which provides asymptotic control of the FDP.

The approach we propose is built upon an underlying procedure that (asymptotically) controls the $k$-FWER for any fixed $k \geq 1$. We then sequentially apply this $k$-FWER procedure for $k = 1, 2, \ldots$ until a stopping rule indicates termination. In the end, we reject all hypotheses that were rejected in the last round of applying the $k$-FWER procedure.

ALGORITHM 4.1 (*Generic method for control of the FDP*).

1. Let $j = 1$ and let $k_1 = 1$.
2. Apply the $k_j$-FWER procedure and denote by $N_j$ the number of hypotheses it rejects.
3.    (a) If $N_j < k_j/\gamma - 1$, stop and reject all hypotheses rejected by the $k_j$-FWER procedure.
   (b) Otherwise, let $j = j + 1$ and then $k_j = k_{j-1} + 1$. Return to Step 2.

Note that the algorithm does not assume anything about the nature of the underlying $k$-FWER procedure. However, in order to reject as many false hypotheses as possible while maintaining (asymptotic) control of the FDP, we suggest to employ a stepwise procedure that accounts for the dependence structure of the test statistics $T_{n,i}$. Algorithm 4.1 is similar to the proposal of [15] for FDP control which is, however, restricted to a multivariate permutation model. The proposal of [15] is heuristic in the sense that they cannot guarantee finite sample or asymptotic control of the FDP even if the permutation hypothesis is valid. However, we will show asymptotic control (and simulations presented later show good finite sample control). The theorem below considers a general bootstrap construction where the individual tests are one-sided and concern univariate parameters $\theta_i(P)$. The bootstrap construction for two-sided tests and the more general subsampling construction can be handled similarly.

THEOREM 4.1. *Consider the setup of Theorem* 3.1. *Fix $P$ satisfying Assumption* B1. *Let $\hat{Q}_n$ be an estimate of $P$ satisfying Assumption* B2. *Employ the step-down procedure of Algorithm* 2.1 *with $\hat{c}_{n,K}(1 - \alpha, k)$ replaced by $b_{n,K}(1 - \alpha, \hat{Q}_n, k)$ as the underlying $k$-FWER procedure. Then the following statements concerning Algorithm* 4.1 *are true:*

   (i) $\limsup_n P\{FDP > \gamma\} \leq \alpha$.
   (ii) *If $P$ is such that $i \notin I(P)$, that is, $H_i$ is false and $\theta_i(P) > 0$, then the probability that the method rejects $H_i$ tends to 1.*

REMARK 4.1. The theorem remains valid if the bootstrap $k$-FWER procedure is based on the operative method of Remark 3.3 or the streamlined

Algorithm 2.2 instead of the generic Algorithm 2.1. But, again, in view of finite sample performance, we suggest the use of the generic Algorithm 2.1 if feasible or at least the use of the operative method.

**5. Comparison with related methods.** We have proposed step-down procedures that control the $k$-FWER and the FDP, with the goal of improving upon methods that do not attempt to incorporate or estimate the dependence structure between the test statistics or $p$-values. An alternative approach toward achieving this goal is given in [34]. We briefly discuss their proposal. (Note that resampling-based procedures of [9] and [34], among others, are implemented in the open source R package `multtest` released as part of the Bioconductor Project; see cran.r-project.org and www.bioconductor.org.)

The approach of [34] begins with an initial procedure that controls the 1-FWER (i.e., the usual FWER) and then rejects in addition the $k-1$ most significant hypotheses not rejected so far. They coin this an *augmentation procedure*, since the 1-FWER rejection set is augmented by the $k-1$ next most significant hypotheses to arrive at the $k$-FWER rejection set. Obviously, if the 1-FWER procedure succeeds in (asymptotically) controlling the 1-FWER, then the augmented procedure provides (asymptotic) control of the $k$-FWER. However, this approach seems suboptimal, because it makes the worst case assumption that, having achieved 1-FWER control, the $k-1$ next most significant hypotheses are all true hypotheses. Moreover, $k-1$ additional hypotheses are always rejected, even if the test statistics or $p$-values to which they correspond are clearly not significant. In addition, the approach really does not fully utilize the weaker measure of error control afforded by using the $k$-FWER with $k > 1$, in that the augmentation method will reject more than $k-1$ hypotheses if and only if the 1-FWER controlling procedure rejects some hypotheses, and this criterion may be too strong to admit any rejections.

Our approach to control the $k$-FWER is based on knowing or estimating the sampling distribution of a suitable $k$-max statistic, that is, the $k$th largest of the $s$ individual (possibly standardized) test statistics. A hypothesis is rejected if its corresponding test statistic is large (relative to the estimated quantiles of the sampling distribution of the $k$-max statistic), unlike the augmentation approach where a hypothesis can be rejected even if its corresponding test statistic is not deemed large by any measure.

To appreciate how the two approaches differ, first consider augmentation based on the Holm procedure, given by (6) with $k = 1$. Other than the additional $k-1$ hypotheses that are rejected after applying Holm, the procedure can only reject a nontrivial number ($k$ or more) if and only if the smallest $p$-value is $\leq \alpha/s$. On the other hand, the generalized Holm procedure starts out with a great advantage; the smallest $p$-value is compared with $k\alpha/s$, a $k$-fold increase. While it is possible for augmentation to reject more hypotheses, it can only reject $k-1$ more than the Holm procedure (and these additional rejections may be suspect because they can correspond to large $p$-values), but the generalized Holm procedure can reject many, many more.

Similar comparisons can be made with augmentation applied to a FWER controlling procedure that attempts to account for the dependence structure (like the ones in this paper with $k=1$). Augmentation might reject $k-1$ more hypotheses than the ones we propose here, but our methods can easily reject many more. Note that, if the test statistics or $p$-values are independent, then augmentation of a bootstrap method that controls the FWER still cannot produce anything much better than the Holm method.

The comparison is similar for the procedures controlling the FDP. Our approach is to sequentially apply a $k$-FWER procedure for $k = 1, 2, \ldots$ until a stopping rule indicates termination. On the other hand, [34] again augment the rejection set of an initial 1-FWER procedure. The idea now is as follows. Let $R$ denote the number of rejections by the 1-FWER procedure. Then reject in addition the $D$ next most significant hypotheses where $D$ is the largest integer which satisfies

$$\frac{D}{D+R} \leq \gamma.$$

Again, if the 1-FWER procedure succeeds in (asymptotically) controlling the 1-FWER, then the augmented procedure provides (asymptotic) control of the FDP. But also again, this approach seems pessimistic in that it makes the worst case assumption that, having achieved 1-FWER control, the $D$ next most significant hypotheses are all true hypotheses.

The next section compares the finite sample performance of the two approaches.

**6. Simulation study.** This section presents a small simulation study in the context of testing population means. We generate random vectors $X_1, \ldots, X_n$ from an $s$-dimensional multivariate normal distribution with mean vector $\theta = (\theta_1, \ldots, \theta_s)$, where $n = 100$ and $s = 50$ or $s = 400$. The null hypotheses are $H_i : \theta_i \leq 0$ and the alternative hypotheses are $H_i : \theta_i > 0$. The test statistics are $T_{n,i} = \sqrt{n}\bar{X}_{i,\cdot}/S_i$, where

$$\bar{X}_{i,\cdot} = \frac{1}{n}\sum_{j=1}^{n} X_{i,j}$$

and

$$S_i^2 = \frac{1}{n-1}\sum_{j=1}^{n}(X_{i,j} - \bar{X}_{i,\cdot})^2.$$

The individual means $\theta_i$ are equal to either 0 or 0.25. The number of means equal to 0.25 is 0, 10, 25 or 50 when $s = 50$ and 0, 100, 200 or 400 when $s = 400$. The covariance matrix is of the common correlation structure $\sigma_{i,i} = 1$ and $\sigma_{i,j} = \rho$ for $i \neq j$. We consider the three values $\rho = 0.0$, 0.5 and 0.8. Other specifications of the covariance matrix do not lead to results that are qualitatively different; see [27].

We include the following multiple testing procedures in the study. The value of $k$ is $k = 3$ when $s = 50$ and $k = 10$ when $s = 400$. The nominal level is $\alpha = 0.05$, unless indicated otherwise.

- (1-Boot) The bootstrap 1-FWER construction of Section 3.1. (This construction is equivalent to the FWER maxT procedure of [9].)
- ($k$-Aug) The $k$-FWER augmentation procedure of [34].
- ($k$-gH) The $k$-FWER generalized Holm procedure described by (6).
- ($k$-Boot) The bootstrap $k$-FWER construction of Section 3.1.
- (Aug$_{0.1}$) The FDP augmentation procedure of [34] with $\gamma = 0.1$.
- (EB$_{0.1}$) The empirical Bayes FDP procedure of [33] with $\gamma = 0.1$.
- (LR$_{0.1}$) The FDP procedure of [17] with $\gamma = 0.1$; see (7).
- (Boot$_{0.1}$) The bootstrap FDP construction of Section 4 with $\gamma = 0.1$.
- (Boot$_{0.1}^{\text{Med}}$) The bootstrap FDP construction of Section 4 with $\gamma = 0.1$ but nominal level $\alpha = 0.5$. Therefore, this procedure asymptotically controls the median FDP to be bounded above by $\gamma = 0.1$.

The augmentation procedures $k$-Aug and Aug$_{0.1}$ are both based on the stepdown 1-Boot construction as the initial 1-FWER controlling procedure. The $k$-Boot procedure is based on the operative method with $N_{\max} = 50$; see Remark 3.3. The estimate $\hat{Q}_n$ employed in the bootstrap is the empirical distribution of the observed data; and for each simulated data set, the same set of $B = 500$ resamples is shared by all bootstrap procedures. The individual $p$-values for $k$-gH and LR$_{0.1}$ are derived from the relation $T_{n,i} \sim t_{n-1}$ under $H_i$.

The performance criteria are (i) the empirical $k$-FWERs and FDPs, compared to the nominal level $\alpha = 0.05$ (or $\alpha = 0.5$ for the method controlling the median FDP); and (ii) the average number of false hypotheses rejected. Since the $k$-Aug procedure rejects the $k - 1$ most significant hypotheses regardless of the data, we also follow this route for the $k$-gH and $k$-Boot procedures to ensure a fair comparison as far as (ii) is concerned (though the differences are really negligible if this route is not followed for the $k$-gH and $k$-Boot procedures). The results are presented in Table 1 for $s = 50$ and in Table 2 for $s = 400$. They can be summarized as follows.

- Almost all methods provide satisfactory finite sample control of their respective $k$-FWER or FDP criteria. In particular, the finite sample control does not appear to deteriorate when the number of hypotheses is increased from $s = 50$ to $s = 400$, while the sample size is kept fixed at $n = 100$.

Table 1

*Empirical FWEs and FDPs expressed as percentages (in the rows "Control") and average number of false hypotheses rejected (in the rows "Rejected") for various methods, with $n = 100$ and $s = 50$. The nominal level is $\alpha = 5\%$, apart from the last column where it is $\alpha = 50\%$. The number of repetitions is 5,000 per scenario and the number of bootstrap resamples is $B = 500$*

| | 1-Boot | 3-Aug | 3-gH | 3-Boot | $Aug_{0.1}$ | $EB_{0.1}$ | $LR_{0.1}$ | $Boot_{0.1}$ | $Boot_{0.1}^{Med}$ |
|---|---|---|---|---|---|---|---|---|---|
| | | | Common correlation: $\rho = 0$ | | | | | | |
| | | | All $\theta_i = 0$ | | | | | | |
| Control | 5.0 | 5.0 | 0.0 | 4.6 | 5.0 | 29.5 | 4.9 | 5.0 | 52.3 |
| Rejected | 0.0 | 0.0 | 0.0 | 0.0 | 0.0 | 0.0 | 0.0 | 0.0 | 0.0 |
| | | | Ten $\theta_i = 0.25$ | | | | | | |
| Control | 4.8 | 0.1 | 0.0 | 3.1 | 4.8 | 24.1 | 4.4 | 4.8 | 47.6 |
| Rejected | 2.6 | 4.5 | 3.9 | 6.3 | 2.6 | 5.0 | 2.6 | 2.6 | 6.3 |
| | | | Twenty-five $\theta_i = 0.25$ | | | | | | |
| Control | 3.3 | 0.0 | 0.0 | 2.2 | 1.9 | 4.5 | 2.1 | 3.0 | 39.1 |
| Rejected | 6.9 | 8.9 | 9.5 | 16.7 | 7.2 | 15.5 | 7.2 | 7.8 | 21.3 |
| | | | All $\theta_i = 0.25$ | | | | | | |
| Control | 0.0 | 0.0 | 0.0 | 0.0 | 0.0 | 0.0 | 0.0 | 0.0 | 0.0 |
| Rejected | 14.9 | 16.9 | 19.2 | 42.3 | 16.2 | 46.8 | 21.3 | 45.3 | 50.0 |
| | | | Common correlation: $\rho = 0.5$ | | | | | | |
| | | | All $\theta_i = 0$ | | | | | | |
| Control | 5.3 | 5.3 | 1.6 | 5.3 | 5.3 | 13.1 | 3.0 | 5.3 | 50.7 |
| Rejected | 0.0 | 0.0 | 0.0 | 0.0 | 0.0 | 0.0 | 0.0 | 0.0 | 0.0 |
| | | | Ten $\theta_i = 0.25$ | | | | | | |
| Control | 5.0 | 2.9 | 1.4 | 4.5 | 4.1 | 9.2 | 2.6 | 4.7 | 49.0 |
| Rejected | 3.4 | 5.2 | 4.3 | 5.6 | 3.4 | 5.0 | 2.7 | 3.4 | 8.3 |
| | | | Twenty-five $\theta_i = 0.25$ | | | | | | |
| Control | 4.3 | 2.0 | 0.1 | 4.4 | 2.8 | 8.4 | 1.6 | 4.5 | 47.2 |
| Rejected | 8.7 | 10.6 | 9.6 | 14.2 | 9.2 | 13.9 | 7.8 | 10.4 | 22.8 |
| | | | All $\theta_i = 0.25$ | | | | | | |
| Control | 0.0 | 0.0 | 0.0 | 0.0 | 0.0 | 0.0 | 0.0 | 0.0 | 0.0 |
| Rejected | 20.2 | 22.0 | 19.2 | 33.0 | 21.5 | 39.3 | 22.5 | 30.6 | 48.9 |
| | | | Common correlation: $\rho = 0.8$ | | | | | | |
| | | | All $\theta_i = 0$ | | | | | | |
| Control | 4.9 | 4.9 | 1.3 | 5.2 | 4.9 | 6.6 | 1.4 | 4.9 | 50.0 |
| Rejected | 0.0 | 0.0 | 0.0 | 0.0 | 0.0 | 0.0 | 0.0 | 0.0 | 0.0 |
| | | | Ten $\theta_i = 0.25$ | | | | | | |
| Control | 5.1 | 5.1 | 1.4 | 4.9 | 4.1 | 6.6 | 1.6 | 4.9 | 49.2 |
| Rejected | 4.8 | 6.3 | 4.6 | 6.4 | 4.9 | 5.1 | 2.7 | 4.9 | 9.3 |
| | | | Twenty-five $\theta_i = 0.25$ | | | | | | |
| Control | 4.6 | 4.6 | 0.1 | 4.5 | 4.5 | 7.4 | 1.6 | 4.5 | 48.0 |
| Rejected | 12.2 | 13.8 | 9.9 | 15.5 | 12.8 | 14.7 | 7.8 | 13.5 | 23.9 |
| | | | All $\theta_i = 0.25$ | | | | | | |
| Control | 0.0 | 0.0 | 0.0 | 0.0 | 0.0 | 0.0 | 0.0 | 0.0 | 0.0 |
| Rejected | 27.1 | 28.4 | 19.5 | 33.7 | 27.9 | 38.0 | 21.4 | 33.1 | 49.0 |

TABLE 2

*Empirical FWEs and FDPs expressed as percentages (in the rows "Control") and average number of false hypotheses rejected (in the rows "Rejected") for various methods, with $n = 100$ and $s = 400$. The nominal level is $\alpha = 5\%$, apart from the last column where it is $\alpha = 50\%$. The number of repetitions is 5,000 when all $\theta_i = 0$ and 2,000 for all other scenarios; and the number of bootstrap resamples is $B = 500$*

| | 1-Boot | 10-Aug | 10-gH | 10-Boot | $\text{Aug}_{0.1}$ | $\text{EB}_{0.1}$ | $\text{LR}_{0.1}$ | $\text{Boot}_{0.1}$ | $\text{Boot}_{0.1}^{\text{Med}}$ |
|---|---|---|---|---|---|---|---|---|---|
| | | | | Common correlation: $\rho = 0$ | | | | | |
| | | | | All $\theta_i = 0$ | | | | | |
| Control | 5.0 | 5.0 | 0.0 | 1.6 | 4.9 | 4.4 | 4.8 | 4.9 | 54.4 |
| Rejected | 0.0 | 0.0 | 0.0 | 0.0 | 0.0 | 0.0 | 0.0 | 0.0 | 0.0 |
| | | | | One hundred $\theta_i = 0.25$ | | | | | |
| Control | 4.3 | 0.0 | 0.0 | 0.5 | 1.0 | 1.7 | 1.5 | 1.7 | 41.0 |
| Rejected | 10.9 | 19.8 | 28.2 | 59.4 | 11.7 | 44.9 | 14.2 | 29.7 | 68.7 |
| | | | | Two hundred $\theta_i = 0.25$ | | | | | |
| Control | 2.7 | 0.0 | 0.0 | 0.4 | 0.0 | 0.1 | 0.0 | 0.4 | 29.9 |
| Rejected | 22.0 | 31.1 | 56.1 | 126.0 | 24.2 | 155.0 | 43.8 | 146.1 | 173.0 |
| | | | | All $\theta_i = 0.25$ | | | | | |
| Control | 0.0 | 0.0 | 0.0 | 0.0 | 0.0 | 0.0 | 0.0 | 0.0 | 0.0 |
| Rejected | 45.4 | 54.4 | 112.4 | 341.1 | 50.5 | 390.0 | 153.7 | 400.0 | 400.0 |
| | | | | Common correlation: $\rho = 0.5$ | | | | | |
| | | | | All $\theta_i = 0$ | | | | | |
| Control | 5.5 | 5.5 | 0.1 | 5.5 | 5.5 | 5.5 | 2.2 | 5.5 | 51.4 |
| Rejected | 0.0 | 0.0 | 0.0 | 0.0 | 0.0 | 0.0 | 0.0 | 0.0 | 0.0 |
| | | | | One hundred $\theta_i = 0.25$ | | | | | |
| Control | 4.9 | 0.4 | 0.5 | 4.4 | 0.5 | 8.0 | 0.7 | 4.2 | 50.5 |
| Rejected | 18.3 | 27.2 | 29.8 | 48.2 | 20.0 | 37.7 | 17.9 | 34.0 | 86.3 |
| | | | | Two hundred $\theta_i = 0.25$ | | | | | |
| Control | 5.0 | 0.4 | 0.5 | 5.1 | 0.3 | 7.8 | 1.1 | 5.0 | 50.2 |
| Rejected | 38.3 | 47.2 | 57.1 | 99.3 | 42.4 | 106.3 | 51.1 | 92.7 | 183.7 |
| | | | | All $\theta_i = 0.25$ | | | | | |
| Control | 0.0 | 0.0 | 0.0 | 0.0 | 0.0 | 0.0 | 0.0 | 0.0 | 0.0 |
| Rejected | 84.3 | 93.3 | 114.0 | 237.7 | 93.0 | 314.9 | 169.8 | 282.7 | 395.5 |
| | | | | Common correlation: $\rho = 0.8$ | | | | | |
| | | | | All $\theta_i = 0$ | | | | | |
| Control | 5.3 | 5.3 | 0.1 | 5.2 | 5.3 | 5.3 | 0.7 | 5.3 | 51.3 |
| Rejected | 0.0 | 0.0 | 0.0 | 0.0 | 0.0 | 0.0 | 0.0 | 0.0 | 0.0 |
| | | | | One hundred $\theta_i = 0.25$ | | | | | |
| Control | 4.9 | 4.1 | 0.5 | 4.5 | 4.9 | 6.2 | 0.7 | 4.5 | 50.8 |
| Rejected | 36.2 | 44.3 | 31.6 | 57.7 | 39.0 | 47.8 | 19.2 | 49.2 | 95.0 |
| | | | | Two hundred $\theta_i = 0.25$ | | | | | |
| Control | 5.4 | 4.2 | 0.1 | 5.4 | 5.4 | 6.6 | 1.2 | 5.4 | 50.5 |
| Rejected | 74.3 | 82.5 | 59.3 | 116.3 | 80.3 | 115.4 | 52.4 | 112.6 | 192.9 |
| | | | | All $\theta_i = 0.25$ | | | | | |
| Control | 0.0 | 0.0 | 0.0 | 0.0 | 0.0 | 0.0 | 0.0 | 0.0 | 0.0 |
| Rejected | 165.7 | 172.8 | 117.0 | 255.1 | 174.4 | 301.9 | 149.5 | 275.3 | 392.8 |

- The exception is $EB_{0.1}$, which can be quite liberal, in particular when $s = 50$ and all null hypotheses are true. As acknowledged to us by the authors of [33], this method is not consistent when all null hypotheses are true and they advocate its use only in settings when false null hypotheses can be anticipated. We provide a brief explanation in Appendix A.

- The relative power of the conservative methods $k$-gH and $LR_{0.1}$ compared to the procedures based on the bootstrap $k$-Aug and $k$-Boot decreases as the common correlation $\rho$ increases.

- Depending on context, $k$-Boot can detect many more false alternatives compared to 1-Boot. The same is not true for $k$-Aug, since, by design, it detects at most $k - 1$ more false hypotheses compared to 1-Boot. So especially when $s$ is large, this approach appears suboptimal. Even the conservative $k$-gH method can be more powerful than the augmentation method for large $s$.

- The comparison is similar for the various FDP procedures. Of all the procedures that provide satisfactory finite sample control, $Boot_{0.1}$ is the most powerful one. $Aug_{0.1}$ becomes uncompetitive when $s$ is large and can even be outperformed by the conservative $LR_{0.1}$ method. Note that $EB_{0.1}$ is often more powerful than $Boot_{0.1}$, but given that its overall finite sample control is not satisfactory, one should be cautious in using this method.

- The power advantage of $k$-Boot and $Boot_{0.1}$ over 1-Boot diminishes as the common correlation $\rho$ increases. (As a result, the same is true for the power advantages of $k$-Boot over $k$-Aug and of $Boot_{0.1}$ over $Aug_{0.1}$, resp.) This is not surprising. Take the extreme case of $\rho = 1$ in our simulation set-up where all nonzero means are equal. In this case 1-Boot rejects either no false hypotheses or all false hypotheses. On the other hand, $k$-Boot rejects either at most $k - 1$ false hypotheses (when the $k - 1$ most significant hypotheses are rejected regardless) or also all false hypotheses. (Note that the 1-max of the "alternative" test statistics will be equal to the $k$-max of the "alternative" test statistics and analogously for the "null" test statistics.) This implies a minimal power gain of $k$-Boot over 1-Boot compared to the case of $\rho = 0$ where the additional number of rejected false hypotheses can far exceed $k - 1$.

The procedure controlling the median FDP (last column) is always the most powerful one. However, it should be understood that it is philosophically different from the other FDP controlling procedures. If $P\{FDP > 0.1\} \leq 0.05$ is achieved, then, in a given application, one can be 95% confident that the realized FDP is at most 0.1. On the other hand, if $P\{FDP > 0.1\} \leq 0.5$ is achieved (i.e., control of the median FDP), then, in a given application, one can only be 50% confident that the realized FDP is at most 0.1. So, loosely speaking, there is a good chance that the realized FDP ends up greater than 0.1, and perhaps by quite a bit. Romano and Wolf [27] examine this issue in more detail by looking at the sampling distribution of the FDP in various scenarios when the median FDP is controlled; see their Figure 1. Depending on the underlying dependence structure, this sampling distribution can exhibit significant variation. As a result, the realized FDP may well be quite above $\gamma = 0.1$.

A similar problem arises in controlling the false discovery rate (FDR), as proposed by [1]. The FDR is the expected value of the FDP. Like the median FDP, it is also a measure of central tendency of the sampling distribution of the FDP. In a given application, the realized FDP can be quite far away from its expected value, the FDR, as made clear in [15].

Finally, some further simulations comparing the augmentation procedures of [34] and the procedures of [17] can be found in [8].

## 7. Concluding remarks.

We have shown how computationally feasible step-down methods can be constructed to control generalized error rates in multiple testing. On the one hand, we have considered the $k$-FWER, which is defined as the probability of making $k$ or more false rejections. This concept would be appropriate when a given number of false rejections can be tolerated. On the other hand, we have also considered the FDP, which is the ratio of false rejections out of the total number of rejections (and defined to be zero when there are no rejections). This concept would be appropriate when a certain proportion of false rejections can be tolerated. Some simulations have shown that these less strict methods can reject many more false hypotheses compared to the traditional FWER control, especially when the number of hypotheses under test is large.

Our step-down methods (asymptotically) account for the dependence structure across test statistics. As a result, they are more powerful than the generalized Holm step-down methods of [14] and [17], which are based on individual $p$-values and designed to handle a "worst case" dependence structure. An alternative approach that also accounts for the dependence structure across test statistics is the augmentation approach of [34]. However, simulations show their methods are noticeably less powerful, especially when the number of hypotheses under test is large. The empirical Bayes method of [33] can sometimes be more powerful than our bootstrap approach for FDP control. However, it also can be quite liberal and it does not offer asymptotic control of the FDP when all null hypotheses are true. Overall, our methods for control of the $k$-FWER and FDP appear competitive with or outperform currently available methods.

## APPENDIX A: PROOFS

PROOF OF THEOREM 2.1. Assume $|I(P)| \geq k$, or there is nothing to prove. Consider the event that at least $k$ true null hypotheses are rejected.

Let $\hat{j}$ be the (random) smallest index $j$ in the algorithm where this occurs, so that $k$-$\max(T_{n,i} : i \in I(P)) > \hat{d}_{n,A_{\hat{j}}}(1 - \alpha, k)$. By definition of $\hat{j}$ (now fixed), $I(P) \subset A_{\hat{j}} \cup I_0$, where $I_0$ is some set of indices satisfying $I_0 \subset R_{\hat{j}}$ and $|I_0| = k - 1$. Let $L$ be any set of indices of false null hypotheses (not necessarily uniquely defined) which satisfy $A_{\hat{j}} \cup I_0 = I(P) \cup L$. Since $\hat{d}_{n,A_{\hat{j}}}(1 - \alpha, k)$ is defined by taking the maximum over sets $I$ of $\hat{c}_{n,K}(1 - \alpha, k)$ with $K = A_{\hat{j}} \cup I$ as $I$ varies over indices satisfying $I \subset R_{\hat{j}}$ and $|I| = k - 1$, it follows that $\hat{d}_{n,A_{\hat{j}}}(1 - \alpha, k) \geq \hat{c}_{n,I(P) \cup L}(1 - \alpha, k)$. By the monotonicity assumption, $\hat{c}_{n,I(P) \cup L}(1 - \alpha, k) \geq \hat{c}_{n,I(P)}(1 - \alpha, k)$. To summarize, the event that at least $k$ true null hypotheses are rejected implies that

$$k\text{-}\max(T_{n,i} : i \in I(P)) > \hat{c}_{n,I(P)}(1 - \alpha, k)$$

and so (i) follows. Part (ii) follows immediately from (i).  □

LEMMA A.1.    *Let $k \leq s$.* (i) *The $k$-$\max$ function is continuous; that is, if $y_n = (y_{n,1}, \ldots, y_{n,s}) \in \mathbb{R}^s$ and $y_n \to y \in \mathbb{R}^s$, then, as $n \to \infty$, $k$-$\max(y_{n,1}, \ldots, y_{n,s}) \to k$-$\max(y_1, \ldots, y_s)$.*

(ii) *If $Y_n \in \mathbb{R}^s$ and $Y_n \xrightarrow{L} Y$, then $k$-$\max(Y_{n,1}, \ldots, Y_{n,s}) \xrightarrow{L} k$-$\max(Y_1, \ldots, Y_s)$.*
(iii) *Furthermore, if each $Y_i$ in (*ii*) has a continuous marginal distribution, then the distribution of $k$-$\max(Y_1, \ldots, Y_s)$ is continuous.*

PROOF.    Part (i) is trivial, and the continuous mapping theorem then implies (ii). To prove (iii), $P\{k$-$\max(Y_1, \ldots, Y_s) = x\} \leq \sum_{i=1}^s P\{Y_i = x\}$.  □

PROOF OF THEOREM 3.1.    To prove (i), by Corollary 3.1 it is sufficient to show that

$$(28) \qquad \limsup_n P\{k\text{-}\max(T_{n,i} : i \in I(P)) > b_{n,I(P)}(1 - \alpha, k, \hat{Q}_n)\} \leq \alpha.$$

Since $\theta_i(P) \leq 0$ for $i \in I(P)$, it follows that

$$k\text{-}\max(T_{n,i} : i \in I(P)) = k\text{-}\max(\tau_n \hat{\theta}_{n,i} : i \in I(P))$$

$$\leq k\text{-}\max(\tau_n[\hat{\theta}_{n,i} - \theta_i(P)] : i \in I(P)).$$

Therefore, the left-hand side of (28) is bounded above by

$$(29) \quad \lim_n P\{k\text{-}\max(\tau_n[\hat{\theta}_{n,i} - \theta_i(P)] : i \in I(P)) > \hat{b}_{n,I(P)}(1 - \alpha, k, \hat{Q}_n)\}.$$

Assumptions B1 and B2 together with the continuous mapping theorem imply that

$$\rho(L_{n,I(P)}(k, P), L_{n,I(P)}(k, \hat{Q}_n)) \xrightarrow{P} 0,$$

for any metric $\rho$ metrizing weak convergence on $\mathbb{R}$. Hence, it follows that (29) is equal to $\alpha$, by an argument very similar to the proof of Theorem 1 of [3].

To prove (ii), assume $\theta_i(P) > 0$. Assumptions B1 and B2 together imply that $b_{n,A_1}(1 - \alpha, k, \hat{Q}_n)$ is stochastically bounded, where $A_1 = \{1, \ldots, s\}$. Furthermore, by the continuous mapping theorem, $\tau_n[\hat{\theta}_{i,n} - \theta_i(P)]$ has a limiting distribution, so $T_{n,i} = \tau_n \hat{\theta}_{i,n} \xrightarrow{P} \infty$. Therefore, with probability tending to one, $T_{n,i} > b_{n,A_1}(1 - \alpha, k, \hat{Q}_n)$, resulting in the rejection of $H_i$ in the first step of Algorithm 2.1. $\quad\square$

PROOF OF THEOREM 3.2. The proof is completely analogous to the proof of Theorem 3.1. The only additional fact needed to prove (iii) is that, when $\theta_i(P) > 0$, $\tau_n \hat{\theta}_{n,i} > 0$ with probability tending to one, and similarly for $\theta_i(P) < 0$. Indeed, Assumption B1(i) implies $\tau_n[\hat{\theta}_{n,i} - \theta_i(P)]$ has a limiting distribution, which implies $\tau_n \hat{\theta}_{n,i} \xrightarrow{P} \infty$ when $\theta_i(P) > 0$, and $\tau_n \hat{\theta}_{n,i} \xrightarrow{P} -\infty$ when $\theta_i(P) < 0$. $\quad\square$

PROOF OF THEOREM 3.3. To prove (i), note that by reasoning similar to before, $\min(T_{n,i} : i \notin I(P)) \xrightarrow{P} \infty$. On the other hand, $\max(T_{n,i} : i \in I(P))$ is either bounded in probability, in case $\theta_i(P) = 0$ for at least one $i \in I(P)$, or $\max(T_{n,i} : i \in I(P)) \xrightarrow{P} -\infty$, in case $\theta_i(P) < 0$ for all $i \in I(P)$. Therefore, the event

$$(30) \qquad \min(T_{n,i} : i \notin I(P)) > \max(T_{n,i} : i \in I(P))$$

has probability tending to 1. But if the event (30) happens, then the rejected true hypotheses (if such exist) will always be the least significant hypotheses among the rejected hypotheses at any stage. This together with the monotonicity of the critical values $b_{n,K}(1 - \alpha, k, \hat{Q}_n)$ allows us to follow asymptotic control of the $k$-FWER from (28) even when Algorithm 2.2 is used. But (28) was already established in the proof of Theorem 3.1.

The proof of (ii) is identical to the proof of (ii) of Theorem 3.1. $\quad\square$

PROOF OF THEOREM 3.4. The proof of (i) is the essential subsampling argument, which derives from (24) being a $U$-statistic; see Theorem 2.6.1 of [21] where one statistic is treated, but the argument is extendable to the simultaneous estimation of the joint distribution. The result (ii) follows as well. To prove (iii), note that by Corollary 3.2 it is sufficient to show that

$$(31) \qquad \limsup_n P\{k\text{-}\max(T_{n,i} : i \in I(P)) > \hat{c}_{n,I(P)}(1 - \alpha, k)\} \leq \alpha.$$

But part (ii) of the theorem implies, for any $\varepsilon > 0$,

$$\hat{c}_{n,I(P)}(1 - \alpha, k) \geq c_{I(P)}(1 - \alpha, k) - \varepsilon \qquad \text{with probability} \to 1.$$

Therefore, using Assumption S, the limit superior of the probability of violation of the $k$-FWER criterion is bounded above, for any $\varepsilon > 0$, by

$$\limsup_n k\text{-FWER}_P \leq P\{k\text{-}\max(T_i, i \in I(P)) > c_{I(P)}(1-\alpha) - \varepsilon\},$$

where $(T_i, i \in I(P))$ denote variables whose joint distribution is $G_{I(P)}(P)$. But letting $\varepsilon \to 0$, the right-hand side of the last expression becomes

$$1 - H_{I(P)}(c_{I(P)}(1-\alpha), P) = 1 - (1-\alpha) = \alpha. \qquad \square$$

PROOF OF THEOREM 4.1. To prove (i), note that by reasoning similar to the proof of part (i) of Theorem 3.3, with probability tending to one, all false hypotheses are rejected before any true hypothesis comes under scrutiny. Therefore, with probability tending to 1, a violation of the FDP criterion occurs if and only if the event

$$(32) \qquad F > \frac{\gamma}{1-\gamma}(s - |I(P)|)$$

occurs, where $F$ is the number of true hypotheses rejected by Algorithm 4.1. Let $F(k)$ denote the number of true hypotheses rejected by the bootstrap $k$-FWER procedure. Furthermore, let $k^*$ denote the smallest integer greater than $(\gamma/(1-\gamma))(s - |I(P)|)$. Assume $|I(P)| \geq k^*$ or there is nothing to prove. By the above argument, we therefore have

$$
\begin{aligned}
\limsup_n P\{\text{FDP} > \gamma\} &= \limsup_n P\{F \geq k^*\} \\
&\leq \limsup_n P\{F(k^*) \geq k^*\} \\
&\leq \alpha \qquad \text{[by part (ii) of Theorem 3.1].}
\end{aligned}
$$
(33)

To see that (33) holds true, note the following two facts. First, the bootstrap $k$-FWER procedure is monotone in $k$: any hypothesis rejected by the $k_1$-FWER procedure will also be rejected by the $k_2$-FWER procedure as long as $k_1 < k_2$. Second, according to step 3(a) of Algorithm 4.1, the algorithm terminates with the application of the $k^*$-FWER procedure, or even before then, if

$$(34) \qquad N_{k^*} < \frac{k^*}{\gamma} - 1.$$

In case all false hypotheses are rejected first, the event (34) happens if and only if

$$(35) \qquad k^* > \frac{\gamma}{1-\gamma}(s - I(P) - [F(k^*) - (k^* - 1)]).$$

By the definition of $k^*$, the inequality (35) will hold as long as $F(k^*) \leq k^* - 1$. Therefore, the event $F(k^*) \leq k^* - 1$ implies that (1) $F(k) \leq k^* - 1$ for any $k < k^*$; and that (2) Algorithm 4.1 terminates with the application of the $k^*$-FWER procedure, or even before then, if all false hypotheses are rejected first (which happens with probability tending to 1). These two facts together demonstrate the validity of (33).

The proof of (i) follows immediately from part (ii) of Theorem 3.1. □

## APPENDIX B

We briefly argue why the method in [33] does not provide even asymptotic control of the FDP when all null hypotheses are true. For this, assume there is one null hypothesis, so $s = 1$ (or $m = 1$ in the notation of [33]); the argument generalizes to arbitrary $s$. Control of the FDP when $s = 1$ reduces to control of the FWER, so the probability of rejecting a true null hypothesis must be bounded above by $\alpha$.

Suppose $X_1, \ldots, X_n$ are i.i.d. $N(\theta, 1)$. Consider testing the null hypothesis $H: \theta = 0$ against $\theta > 0$. Let $T_n = T_{n,1} = n^{-1/2} \sum_{i=1}^n X_i$. Let $\Phi(\cdot)$ denote the c.d.f. of the standard normal distribution and $\phi(\cdot)$ its density.

Under $H$, $T_n \sim N(0, 1)$ and so $f_0$ (in the notation of [33]) is $\phi$. The algorithm of [33] simplifies to the following:

1. If $T_n \leq 0$, let $\pi = 1$; otherwise, let $\pi = \phi(T_n)/\phi(0)$.
2. Determine $c$ as the solution to $\pi(1 - \Phi(c)) = \alpha$. That is, $c = \Phi^{-1}(1 - \alpha/\pi)$, with $\Phi^{-1}(\lambda)$ defined as $-\infty$ for $\lambda \leq 0$.
3. Reject $H$ if $T_n > c$.

Now assume $\theta = 0$. Since $\Phi^{-1}(1 - \alpha/\pi) < \Phi^{-1}(1 - \alpha)$ with positive probability, $H$ is rejected with probability greater than $\alpha$. Moreover, $\pi$ does not even converge to 1 in probability [since $T_n$ has a nondegenerate $N(0, 1)$ distribution for every $n$, and asymptotically in typical nonparametric problems]. As an example, for $\alpha = 0.05$, a numerical simulation based on 100,000 repetitions results in a rejection probability of 0.107.

**Acknowledgments.** We would like to thank the Editor, the Associate Editor, and two anonymous referees for their careful reviews of this manuscript which have led to an improved presentation.

| DEPARTMENT OF STATISTICS | INSTITUTE FOR EMPIRICAL RESEARCH IN ECONOMICS |
|---|---|
| STANFORD UNIVERSITY | UNIVERSITY OF ZURICH |
| STANFORD, CALIFORNIA 94305-4065 | CH-8006 ZÜRICH |
| USA | SWITZERLAND |
| E-MAIL: romano@stanford.edu | E-MAIL: mwolf@iew.unizh.ch |